\documentclass[a4paper,11pt,twoside,reqno]{amsart}

\usepackage[utf8]{inputenc}
\usepackage[plainpages=false,pdfpagelabels=true]{hyperref}
\usepackage{amssymb,amsthm,wasysym}
\usepackage{dsfont}
\usepackage{slashed}

\numberwithin{equation}{section}
\allowdisplaybreaks[1]

\newtheorem{Satz}{Theorem}[section]
\newtheorem{Prop}[Satz]{Proposition}
\newtheorem{Lem}[Satz]{Lemma}

\newtheorem{Cor}[Satz]{Corollary}
\newcommand{\cR}{{\mathcal R}}

\theoremstyle{definition}
\newtheorem{Dfn}[Satz]{Definition}
\newtheorem{Bem}[Satz]{Remark}

\renewcommand{\epsilon}{\varepsilon}

\newcommand{\R}{\ensuremath{\mathbb{R}}}
\newcommand{\C}{\ensuremath{\mathbb{C}}}
\newcommand{\Z}{\ensuremath{\mathbb{Z}}}

\newcommand{\D}{\slashed{D}}
\newcommand{\p}{\slashed{\partial}}
\newcommand{\sff}{\mathrm{I\!I}}
\newcommand{\Dtor}{\D^{\scriptscriptstyle Tor}}
\newcommand{\nablator}{\nabla^{\scriptscriptstyle Tor}}
\newcommand{\nablattor}{\tilde{\nabla}^{\scriptscriptstyle Tor}}

\title{Dirac-harmonic maps with torsion}
\author{Volker Branding}
\date{\today}
\address{TU Wien\\
Institut für diskrete Mathematik und Geometrie\\
Wiedner Hauptstraße 8–10, A-1040 Wien}
\email[]{volker@geometrie.tuwien.ac.at}
\begin{document}

\begin{abstract}
We study Dirac-harmonic maps from surfaces to manifolds with torsion,
which is motivated from the superstring action considered in theoretical physics.
We discuss analytic and geometric properties of such maps and
outline an existence result for uncoupled solutions.
\end{abstract} 

\maketitle
  
\section{Introduction and Results}
\emph{Dirac-harmonic maps} arise as the mathematical version of the simplest supersymmetric non-linear sigma model studied 
in quantum field theory. They are critical points of an energy functional that couples the equation for harmonic maps
to so-called vector spinors \cite{MR2262709}. If the domain is two-dimensional, Dirac-harmonic maps belong to the class
of conformally invariant variational problems.

Many results for Dirac-harmonic maps have already been obtained. 
This includes the regularity of solutions \cite{MR2176464}, \cite{MR2506243}, \cite{MR2544729} and the energy identity \cite{MR2176464}.
In addition, an existence result for uncoupled solutions \cite{ammannginoux}, for the boundary value problem \cite{springerlink:10.1007/s00526-012-0512-5}, \cite{BVP} 
and for a non-linear version of Dirac-geodesics \cite{springerlink:10.1007/s00526-011-0404-0} have been established.
A heat flow approach for Dirac-harmonic maps has been studied in \cite{phd}, see also \cite{regularized}.

However, in quantum field theory more complicated models are studied.
Taking into account an additional curvature term in the energy functional one is led to \emph{Dirac-harmonic maps with curvature term}, see \cite{MR2370260}.
From an analytical point of view the latter are more difficult and not much is known about solutions of these equations.
Dirac-harmonic maps coupled to a two-form potential, called \emph{Magnetic Dirac-harmonic maps}, are studied in \cite{magnetic}.

The full \((1,1)\) supersymmetric nonlinear \(\sigma\)-model considered in theoretical physics
involves additional terms that are not captured by the previous analysis. 
Some of these additional terms can be interpreted as considering both Dirac-harmonic maps and Dirac-harmonic maps with curvature term
into manifolds having a connection with torsion.

In this note we want to extend the framework of Dirac-harmonic maps to target spaces with torsion.
It turns out that most of the known results for Dirac-harmonic maps still hold,
in particular the regularity of weak solutions and the removable singularity theorem.
Moreover, we outline an approach to the existence question for Dirac-harmonic maps with torsion using index theory.

This paper is organized as follows. In the second section we provide some background material
on the superspace formalism used in theoretical physics and briefly review orthogonal connections with torsion dating back to Cartan.
Section three then introduces Dirac-harmonic maps with torsion and afterwards we discuss geometric (Section 4) and analytic aspects (Section 5)
of these maps. In the last section we comment on Dirac-harmonic maps with curvature term to target manifolds with torsion.
\par\medskip
\textbf{Acknowledgements:}
The author would like to thank Christoph Stephan and Florian Hanisch for several discussions about torsion and supergeometry.

\section{Some Background Material}
\subsection{The (1,1) supersymmetric nonlinear \texorpdfstring{\(\sigma\)}{}-model in superspace}
In this section we want to give a short overview on how physicists formulate supersymmetric sigma models as field theories in superspace.
For a detailed discussion we refer to the books \cite{MR1707282} and \cite{MR1701598},
for more specific details of the \((1,1)\) supersymmetric \(\sigma\)-model
one may consult \cite{MR2151030}, p.106, \cite{Callan:1989nz}, Chapter 5 and references therein.

In two-dimensional superspace we have the usual commuting coordinates \(\xi_+,\xi_-\) and in addition anti-commuting coordinates \(\theta_+,\theta_-\).
The central objects are the so-called superfields \(\Phi(\xi,\theta_+,\theta_-)\), whose components are given in terms of local coordinates by
\begin{equation}
\Phi^j(\xi,\theta_+,\theta_-)=\phi^j(\xi)-i\theta_-\psi_+^j(\xi)+i\theta_+\psi_-^j(\xi).
\end{equation}
Here, \(\phi(\xi)\) denotes a usual map and \(\psi_+(\xi),\psi_-(\xi)\) are certain spinors taking values in a Grassmann algebra. 
We have neglected any auxiliary fields. 
To obtain an action functional, we need the supercovariant derivatives
\[
D_\pm=i\frac{\partial}{\partial\theta_\mp}+\theta_\mp\frac{\partial}{\partial\xi_\pm}.
\]
Using the metric \(g\) on the target manifold we obtain a conformal invariant action for the \((1,1)\) supersymmetric \(\sigma\)-model in superspace 
by setting
\begin{equation}
E^{(1,1)}_{SYM}(\Phi)=\frac{1}{2}\int g(D_+\Phi,D_-\Phi) d^2\xi d\theta_+d\theta_-.
\end{equation}

Expanding the superfield and interpreting the terms from a geometric point of view yields
(the precise definition of all terms is given in Section 3)
\begin{equation}
\label{energy-symmetric}
E(\phi,\psi)_{SYM}=\frac{1}{2}\int_M|d\phi|^2+\langle\psi,\D\psi\rangle+\frac{1}{6}\langle R^N(\psi,\psi)\psi,\psi\rangle.
\end{equation}
The first two terms in the functional give rise to the energy for \emph{Dirac-harmonic maps},
including also the third terms leads to \emph{Dirac-harmonic maps with curvature term}.

But there is another way to write down a conformal invariant action for a supersymmetric nonlinear \(\sigma\)-model in superspace
using a two-form \(B\) on the target manifold. More precisely, one studies the action
\begin{equation}
E^{(1,1)}_{ASYM}(\Phi)=\frac{1}{2}\int B(D_+\Phi,D_-\Phi) d^2\xi d\theta_+d\theta_-.
\end{equation}
Again, we may expand this action in terms of ordinary fields, which gives
\begin{align}
\label{easym}
E(\phi,\psi)_{ASYM}=\frac{1}{2}\int_M &\phi^{-1}B+C(e_\alpha\cdot\psi,\psi,d\phi(e_\alpha))+H(\psi,\psi,\psi,\psi)
\end{align}
with a two-form \(B\), a three-form \(C\) and some quantity \(H\).
The geometric version of the full \((1,1)\) supersymmetric nonlinear \(\sigma\)-model is then governed by the action
\[
E(\phi,\psi)=E(\phi,\psi)_{SYM}+E(\phi,\psi)_{ASYM}.
\]
We want to analyze this action from the point of view of differential geometry.

\subsection{A Shortcut to Torsion}
Orthogonal connections with torsion have already been classified by Cartan, see \cite{MR1509253}, \cite{MR1509255} and \cite{MR1509263}.
However, here we mostly follow the presentation from \cite{MR2873850}, Section 2.

Consider a manifold \(N\) with a Riemannian metric \(g\). By \(\nabla^{\scriptscriptstyle LC}\) we denote
the Levi-Civita connection. For any affine connection there exists a \((2,1)\)-tensor field \(A\)
such that
\begin{equation}
\label{nabla-torsion}
\nablator_XY=\nabla^{\scriptscriptstyle LC}_XY+A(X,Y)
\end{equation}
for all vector fields \(X,Y\in\Gamma(TN)\).
We demand that the connection is \emph{orthogonal}, that is for all vector fields \(X,Y,Z\) one has
\begin{equation}
\label{nabla-metric}
\partial_X\langle Y,Z\rangle =\langle\nabla_XY,Z\rangle+\langle Y,\nabla_XZ\rangle,
\end{equation}
where \(\langle\cdot,\cdot\rangle\) denotes the scalar product of the metric \(g\).
Combing \eqref{nabla-torsion} and \eqref{nabla-metric} we follow that the endomorphism \(A(X,\cdot)\) is skew-adjoint, that is
\begin{equation}
\label{torsion-skewadjointess}
\langle A(X,Y),Z\rangle=-\langle Y,A(X,Z)\rangle.
\end{equation}
The curvature tensors of \(\nabla^{\scriptscriptstyle LC}\) and \(\nablator\) satisfy the following relation
\begin{align}
\label{curvature-tensors}
R^{\scriptscriptstyle Tor}(X,Y)Z=&R^{\scriptscriptstyle LC}(X,Y)Z
+(\nabla^{\scriptscriptstyle LC}_XA)(Y,Z)-(\nabla^{\scriptscriptstyle LC}_YA)(X,Z) \\
\nonumber &+A(X,A(Y,Z))-A(Y,A(X,Z)).
\end{align}	
Regarding the symmetries of the curvature tensor of an orthogonal connection with torsion, we have
(with \(X,Y,Z,W\in\Gamma(TN)\))
\begin{align*}
\langle R^{\scriptscriptstyle Tor}(X,Y)Z,W\rangle=-\langle R^{\scriptscriptstyle Tor}(Y,X)Z,W\rangle, \\
\langle R^{\scriptscriptstyle Tor}(X,Y)Z,W\rangle=-\langle R^{\scriptscriptstyle Tor}(X,Y)W,Z\rangle.
\end{align*}	
However, in general the curvature tensor is not symmetric under swapping the first two entries with the last two entries,
see \cite{MR2322400}, Remark 2.3.

Any torsion tensor \(A\) induces a \((3,0)\) tensor by setting
\[
A_{XYZ}=\langle A(X,Y),Z\rangle.
\]
We define the space of all possible torsion tensors on \(T_pN\) by
\[
\mathcal{T}(T_pN)=\big\{A\in\otimes^3T^*_pN\mid A_{XYZ}=-A_{XZY}~ X,Y,Z\in T_pN\big\}.
\]
For \(A\in\mathcal{T}(T_pM)\) and \(Z\in T_pM\) one sets
\[
c_{12}(A)(Z)=A_{\partial_{y^i}\partial_{y^i}Z},
\]
where \(\partial_{y^i}\) is a local basis of \(TN\) and we sum over \(i\).
The following classification result is due to Cartan:
\begin{Satz}
\label{torsion-classification}
Assume that \(\dim N\geq 3\). Then the space \(\mathcal{T}(T_pN)\) has the following irreducible
decomposition
\[
\mathcal{T}(T_pN)=\mathcal{T}_1(T_pN)\oplus\mathcal{T}_2(T_pN)\oplus\mathcal{T}_3(T_pN),
\]
which is orthogonal with respect to \(\langle\cdot,\cdot\rangle\) and is explicitly given by
\begin{align*}
\mathcal{T}_1(T_pN)=& \{A\in\mathcal{T}(T_pN)\mid\exists V\textrm{ s.t. } A_{XYZ}=\langle X,Y\rangle\langle V,Z\rangle-\langle X,Z\rangle\langle V,Y\rangle\}, \\
\mathcal{T}_2(T_pN)=& \{A\in\mathcal{T}(T_pN)\mid A_{XYZ}=-A_{YXZ}~~\forall X,Y,Z\}, \\
\mathcal{T}_3(T_pN)=& \{A\in\mathcal{T}(T_pN)\mid A_{XYZ}+A_{YZX}+A_{ZXY}=0,~~ c_{12}(A)(Z)=0\}. 
\end{align*}
Moreover, for \(\dim N=2\) we have
\[
\mathcal{T}(T_pN)=\mathcal{T}_1(T_pN).
\]
\end{Satz}
A proof of the above Theorem can be found in \cite{MR712664}, Theorem 3.1.

We call the torsion of a connection, whose torsion tensor is contained in \(\mathcal{T}_1(T_pN)\) \emph{vectorial},
with torsion tensor in \(\mathcal{T}_2(T_pN)=\Lambda^3T^\ast_pN\) \emph{totally anti-symmetric}
and with torsion tensor in \(\mathcal{T}_3(T_pN)\) \emph{of Cartan type}.

In terms of local coordinates we have from \eqref{torsion-skewadjointess}
\begin{equation}
A_{ijk}=-A_{ikj} 
\end{equation}
and from \eqref{curvature-tensors}
\begin{equation}
R^{\scriptscriptstyle Tor}_{ijkl}=R^{\scriptscriptstyle LC}_{ijkl}+\nabla_i A_{jkl}-\nabla_j A_{ikl}
+A_{irl}A_{jk}^{\hspace{0.3cm}r}-A_{jrl}A_{ik}^{\hspace{0.3cm}r}.
\end{equation}

For more details on the geometric/physical interpretation of manifolds with torsion we refer to the lecture notes \cite{MR2322400}
and the survey article \cite{MR1883845}.

\section{Dirac-harmonic maps with torsion}
Let us now describe the geometric framework for Dirac-harmonic maps with torsion in detail.
We assume that \((M,h)\) is a closed Riemannian spin surface with spinor bundle \(\Sigma M\)
and \((N,g)\) is a compact Riemannian manifold. Let \(\phi\colon M\to N\) be a map.
Together with the pull-back bundle \(\phi^{-1}TN\) we may consider the twisted bundle
\(\Sigma M\otimes\phi^{-1}TN\). Sections in this bundle are called \emph{vector spinors},
in terms of local coordinates \(y^i\) on \(N\) they can be expressed as
\begin{equation}
\psi=\psi^i\otimes\frac{\partial}{\partial y^i}(\phi(x)).
\end{equation}
Note that these spinors do not take values in a Grassmann algebra.
We are using the Einstein summation convention, that is, we sum over repeated indices.
Indices on \(N\) will be denoted by Latin letters, whereas indices on \(M\) are labeled by Greek letters.
On \(\Sigma M\otimes\phi^{-1}TN\) we have a connection that is induced 
from the connections on \(\Sigma M\) and \(\phi^{-1}TN\),
we will denote this connection by \(\tilde{\nabla}\).
We will mostly be interested in connections with torsion on the target manifold \(N\),
in this case we will write \(\nablattor\) and then we have the following decomposition
\begin{equation}
\nablattor=\nabla^{\Sigma M}\otimes\mathds{1}^{\phi^{-1}TN} +\mathds{1}^{\Sigma M}\otimes\nabla^{\phi^{-1}TN}+\mathds{1}^{\Sigma M}\otimes A(\cdot,\cdot).
\end{equation}
On the spinor bundle \(\Sigma M\) we have the Clifford multiplication with tangent vectors, which is skew-symmetric
\[
\langle X\cdot\psi,\chi\rangle_{\Sigma M}=-\langle \psi,X\cdot\chi\rangle_{\Sigma M}
\]
for all \(X\in TM\) and \(\psi,\chi\in\Gamma(\Sigma M)\).

We now consider the twisted Dirac operator on \(\Sigma M\otimes\phi^{-1}TN\), namely
\[
\Dtor:=e_\alpha\cdot\nablattor_{e_\alpha}=\D+e_\alpha\cdot A(d\phi(e_\alpha),\cdot),
\]
where \(\{e_\alpha\}\) denotes a local orthonormal basis of \(TM\).
This operator is elliptic and self-adjoint with respect to the \(L^2\) norm, since the connection on \(\phi^{-1}TN\) is metric.
In terms of local coordinates we may express it as
\[
\Dtor\psi=\p\psi^i\otimes\frac{\partial}{\partial y^i}+e_\alpha\cdot\psi^k\otimes\Gamma^i_{jk}\frac{\partial\phi^j}{\partial x_\alpha}\frac{\partial}{\partial y^i}
+e_\alpha\cdot\psi^k\otimes A_{jk}^{\hspace{0.3cm}i}\frac{\partial\phi^j}{\partial x_\alpha}\frac{\partial}{\partial y^i}
\]
with the Christoffel symbols \(\Gamma^i_{jk}\) and the torsion coefficients \(A_{jk}^{\hspace{0.3cm}i}\) on \(N\).
Moreover, \(\p\) denotes the usual Dirac operator acting on sections of \(\Sigma M\).

We may now study the energy functional
\begin{align}
\label{energy-torsion}
E_{\scriptscriptstyle{Tor}}(\phi,\psi)&=\frac{1}{2}\int_M|d\phi|^2+\langle\psi,\Dtor\psi\rangle \\
\nonumber &=\frac{1}{2}\int_M|d\phi|^2+\langle\psi,\D\psi\rangle +
\langle\psi,A(d\phi(e_\alpha),e_\alpha\cdot\psi)\rangle,
\end{align}
which is part of the full \((1,1)\) supersymmetric non-linear sigma model \eqref{easym}
as described in the introduction.

\begin{Bem}
The energy functional \eqref{energy-torsion} is real-valued.
One the one hand this follows from the fact that the operator \(\Dtor\) is elliptic and self-adjoint,
on the other hand we note
\begin{align*}
\overline{\langle\psi,e_\alpha\cdot A(d\phi(e_\alpha),\psi)\rangle}
&=\langle e_\alpha\cdot A(d\phi(e_\alpha),\psi),\psi\rangle
=-\langle A(d\phi(e_\alpha),\psi),e_\alpha\cdot\psi\rangle \\
&=\langle\psi,e_\alpha\cdot A(d\phi(e_\alpha),\psi)\rangle,
\end{align*}
where we used the skew-symmetry of the Clifford multiplication and
the skew-adjointness of the endomorphism \(A\).
\end{Bem}

\begin{Bem}
We only consider a connection with torsion on the target manifold \(N\).
Of course, we could also consider a connection with torsion on the domain \(M\). 
It is known that in this case the Dirac operator is still self-adjoint if
\[
A\in\mathcal{T}_2(TM)\oplus\mathcal{T}_3(TM),
\]
see \cite{MR547866}, Satz 2 and also \cite{MR2873850}, Cor. 4.6.
However, by Theorem \ref{torsion-classification} we know that in the case of a two-dimensional domain \(M\)
only the vectorial torsion contributes. Hence, we would get a twisted Dirac operator
which is no longer self-adjoint. 
\end{Bem}

As for Dirac-harmonic maps we can compute the critical points of \eqref{energy-torsion}:

\begin{Prop}
\label{euler-lagrange-dh}
The critical points of the functional \eqref{energy-torsion} are given by
\begin{align}
\label{euler-lagrange-phi-torsion}\tau(\phi)&=\cR(\phi,\psi)+F^{\scriptscriptstyle Tor}(\phi,\psi),\\
\label{euler-lagrange-psi-torsion}\Dtor\psi&=0
\end{align}
with the curvature term
\begin{align*}
\cR^{\scriptscriptstyle}(\phi,\psi)&=\frac{1}{2}R^N_{\scriptscriptstyle}(e_\alpha\cdot\psi,\psi)d\phi(e_\alpha)
\end{align*}
and the torsion term \(F^{\scriptscriptstyle Tor}(\phi,\psi)\in\Gamma(\phi^{-1}TN)\) defined by \eqref{definition-ftor}.
\end{Prop}

\begin{proof}
We choose a local orthonormal basis \(\{e_\alpha\}\) on \(M\) such that
\([e_\alpha,\partial_t]=0\) and also \(\nabla_{\partial_t}e_\alpha=0\) at a considered point.
Consider a smooth variation of the pair \((\phi,\psi)\) satisfying \((\frac{\partial\phi_t}{\partial t},\frac{\nablattor\psi_t}{\partial t})\big|_{t=0}=(\eta,\xi)\).
Using the skew-adjointness of the endomorphism \(A\), we find
\begin{align*}
\frac{\partial}{\partial t}\big|_{t=0}\frac{1}{2}\int_M|d\phi_t|^2&=\int_M\langle\frac{\nabla^{\scriptscriptstyle LC}}{\partial t}d\phi_t(e_\alpha),d\phi_t(e_\alpha)\rangle\big|_{t=0} \\
&=\int_M\langle\nabla^{\scriptscriptstyle LC}_{e_\alpha}\frac{\partial\phi_t}{\partial t},d\phi_t(e_\alpha)\rangle \big|_{t=0}\\
&=-\int_M\langle\frac{\partial\phi_t}{\partial t},\nabla^{\scriptscriptstyle LC}_{e_\alpha} d\phi_t(e_\alpha)\rangle\big|_{t=0}
=-\int_M\langle\tau(\phi),\eta\rangle.
\end{align*}
Moreover, we calculate
\begin{align*}
\frac{\partial}{\partial t}\frac{1}{2}\int_M&\langle\psi_t,\Dtor\psi_t\rangle \\
&=\frac{1}{2}\int_M\langle\frac{\tilde{\nabla}^{\scriptscriptstyle{LC}}\psi_t}{\partial t},\Dtor\psi_t\rangle
+\langle\psi_t,\frac{\tilde{\nabla}^{\scriptscriptstyle{LC}}}{\partial t}\Dtor\psi_t\rangle \\
&=\frac{1}{2}\int_M\langle\frac{\nablattor\psi_t}{\partial t},\Dtor\psi_t\rangle
+\langle\psi_t,\frac{\nablattor}{\partial t}\Dtor\psi_t\rangle \\
&=\int_M\operatorname{Re}\langle\frac{\nablattor\psi_t}{\partial t},\Dtor\psi_t\rangle
+\frac{1}{2}\langle\psi_t,e_\alpha\cdot R^N_{\scriptscriptstyle Tor}(d\phi_t(\partial_t),d\phi_t(e_\alpha))\psi_t\rangle.
\end{align*}	
Expanding the curvature tensor using \eqref{curvature-tensors}, we find
\begin{align*}
&\langle\psi_t,e_\alpha\cdot R^N_{\scriptscriptstyle Tor}(d\phi_t(\partial_t),d\phi_t(e_\alpha))\psi_t\rangle=\langle\psi_t,e_\alpha\cdot R^N(d\phi_t(\partial_t),d\phi_t(e_\alpha))\psi_t\rangle \\
&+\langle \psi_t,e_\alpha\cdot (\nabla_{d\phi_t(\partial_t)}A)(d\phi_t(e_\alpha),\psi_t)\rangle-
\langle \psi_t,e_\alpha\cdot (\nabla_{d\phi_t(e_\alpha)}A)(d\phi_t(\partial_t),\psi_t)\rangle \\
&+\langle\psi_t,e_\alpha\cdot A(d\phi_t(\partial_t),A(d\phi_t(e_\alpha),\psi_t))\rangle
-\langle\psi_t,e_\alpha\cdot A(d\phi_t(e_\alpha),A(d\phi_t(\partial_t),\psi_t))\rangle.	
\end{align*}
Using the symmetries of the curvature tensor without torsion, we get
\[
\langle\psi_t,e_\alpha\cdot R^N(d\phi_t(\partial_t),d\phi_t(e_\alpha))\psi_t\rangle\big|_{t=0}=\langle R^N_{\scriptscriptstyle}(e_\alpha\cdot\psi,\psi)d\phi(e_\alpha),\eta\rangle.
\]	
For the rest of the terms we define  \(F^{\scriptscriptstyle Tor}(\phi,\psi)\in\Gamma(\phi^{-1}TN)\) by
\begin{align}
\label{definition-ftor}
\langle F^{\scriptscriptstyle Tor}&(\phi,\psi),\eta\rangle:=
\nonumber\frac{1}{2}\big(\langle\psi,e_\alpha\cdot(\nabla_\eta A)(d\phi(e_\alpha),\psi)\rangle
-\langle\psi,e_\alpha\cdot(\nabla_{d\phi(e_\alpha)} A)(\eta,\psi)\rangle \\
&+\langle\psi,e_\alpha\cdot A(\eta,A(d\phi(e_\alpha),\psi))\rangle
-\langle\psi,e_\alpha\cdot A(d\phi(e_\alpha),A(\eta,\psi))\rangle\big)
\end{align}
and evaluating at \(t=0\) 
\begin{equation}
\label{variation-energy}
\frac{d}{dt}E_{\scriptscriptstyle Tor}(\phi_t,\psi_t)\big|_{t=0}=\int_M\operatorname{Re}\langle\xi,\Dtor\psi_t\rangle
+\langle\eta,-\tau(\phi)+\cR(\phi,\psi)+F^{\scriptscriptstyle Tor}(\phi,\psi)\rangle
\end{equation}
gives the result.
\end{proof}

We call solutions \((\phi,\psi)\) of the system \eqref{euler-lagrange-phi-torsion} and \eqref{euler-lagrange-psi-torsion}
\emph{Dirac-harmonic maps with torsion}.

Expanding the connection on \(N\), we find 
\begin{align}
\tau(\phi)=&\frac{1}{2}R^N(e_\alpha\cdot\psi,\psi)d\phi(e_\alpha)+F^{\scriptscriptstyle Tor}(\phi,\psi),\\
\D\psi=&-A(d\phi(e_\alpha),e_\alpha\cdot\psi).
\end{align}
For a general torsion tensor \(A\) the expression \(F^{\scriptscriptstyle Tor}(\phi,\psi)\) cannot be brought into a ``nicer'' form.
However, for vectorial torsion we find
\begin{align*}
F^{\scriptscriptstyle Tor}(\phi,\psi)=&\langle V,d\phi(e_\alpha)\rangle\langle V,e_\alpha\cdot\psi\rangle\psi
-|V|^2\langle d\phi(e_\alpha),e_\alpha\cdot\psi\rangle\psi \\
&+\langle V,\psi\rangle\langle d\phi(e_\alpha),e_\alpha\cdot\psi\rangle V
-\langle d\phi(e_\alpha),e_\alpha\cdot\psi\rangle\langle\psi,(\nabla V)^\sharp\rangle \\
&-\langle\nabla_{d\phi(e_\alpha)}V,e_\alpha\cdot\psi\rangle\psi,
\end{align*}
where \(V\) is a vector field on \(N\).

In terms of local coordinates \(x_\alpha\) on \(M\), the equations for Dirac-harmonic maps with torsion \eqref{euler-lagrange-phi-torsion} and \eqref{euler-lagrange-psi-torsion} acquire the form
\begin{align*}
\tau^m(\phi)=
&\frac{1}{2}R^m_{~lij}\langle\psi^i, e_\alpha\cdot\psi^j\rangle_{\Sigma M}\frac{\partial\phi^l}{\partial x_\alpha} \\
&+\frac{1}{2}\big(\nabla^mA_{lji}-\nabla_lA^m_{~~~ji}+A^m_{~~~~ri}A_{lj}^{\hspace{0.25cm}r}-A_{lri}A^{m\hspace{0.1cm}r}_{~~~~j}\big)\langle\psi^i, e_\alpha\cdot\psi^j\rangle_{\Sigma M}\frac{\partial\phi^l}{\partial x_\alpha},\\
\p\psi^i=&-(A_{jk}^{\hspace{0.3cm} i}+\Gamma^i_{jk})e_\alpha\cdot\psi^j\frac{\partial\phi^k}{\partial x_\alpha}.
\end{align*}

\begin{Bem}
We do not get a torsion contribution for the tension field \(\tau(\phi)\) when starting
from a variational principle. However, if we just take the harmonic map equation and change to a connection
with torsion, then we do get a contribution. In this case the torsion piece
in the tension field vanishes for totally antisymmetric torsion due to symmetry reasons. 
This is the reason why physical models usually consider only skew-symmetric torsion.
\end{Bem}

We call a solution of the Euler-Lagrange equations \eqref{euler-lagrange-phi-torsion} and \eqref{euler-lagrange-psi-torsion}
uncoupled, if \(\phi\) is a harmonic map.
\par\medskip
Using tools from index theory, a general existence result for uncoupled Dirac-harmonic maps
could be derived in \cite{ammannginoux}. Since the index of the twisted Dirac-operator does not 
change when considering a connection with torsion on \(\phi^{-1}TN\) the arguments from \cite{ammannginoux} can also be applied 
in our case. Thus, let us briefly recall the following facts:

Let \((M,h)\) be a closed Riemannian spin manifold of dimension \(m\) with spin structure \(\sigma\) 
and let \(E\to M\) be a vector bundle with \emph{metric} connection. 
The twisted Dirac-operator \(\D^E\colon\Gamma(\Sigma M\otimes E)\to\Gamma(\Sigma M\otimes E)\) 
has an index \(\alpha(M,\sigma,E)\in\operatorname{KO}_m(pt)\)(\cite{MR1031992}, p.141, p.151), where
\begin{align*}
\operatorname{KO}_m(pt)\cong
\begin{cases} 
\Z &\mbox{if } m= 0(4) \\
\Z_2 &\mbox{if } m= 1,2(8) \\
0 &\mbox{otherwise}.
\end{cases}
\end{align*}
On the other hand, the index \(\alpha(M,\sigma,E)\) can be calculated from \(\ker\D^E\)
using \cite{MR1031992}, Thm. 7.13 (with \(\operatorname{ch(E)}\) being the Chern character of the bundle \(E\)):
\begin{align*}
\alpha(M,\sigma,E)=
\begin{cases} 
\{\operatorname{ch(E)}\cdot\widehat A(TM)\}[M]&\mbox{if } m= 0(8) \\
[\dim_\C(\ker(\D^E))]_{\Z_2} &\mbox{if } m= 1(8) \\
[\frac{\dim_\C(\ker(\D^E))}{2}]_{\Z_2}&\mbox{if } m= 2(8) \\
\frac{1}{2}\{\operatorname{ch(E)}\cdot\widehat A(TM)\}[M] &\mbox{if } m= 4(8) \\
\end{cases}
\end{align*}
These statements still hold in our case since we are assuming that we have a metric
connection on \(E=\phi^{-1}TN\). 
From the variational formula \eqref{variation-energy} it can be deduced that in order to obtain an
existence result we have to do the following:
For a given harmonic map \(\phi_0\) and \(\psi_0\in\ker(\Dtor_0)\), where \(\Dtor_0\psi\in\Gamma(\Sigma M\otimes\phi^{-1}_0TN)\),
we have to construct for any smooth variation \(\phi_t\) of \(\phi_0\) a smooth variation of \(\psi_t\) satisfying 
\(\frac{d}{dt}\int_M\langle\psi_t,\Dtor_t\psi_t\rangle\big|_{t=0}=0\).
This is the same argument as Cor. 5.2 in \cite{ammannginoux}.
Note that \(\D\) and \(\Dtor\) have the same principal symbol and the same index.
Hence, this smooth variation can be constructed by assuming that the index \(\alpha(M,\sigma,E)\) is non-trivial, see \cite{ammannginoux}, Prop. 8.2 and Section 9.

\section{Geometric Aspects of solutions}
In this section we analyze some geometric properties of Dirac-harmonic maps with torsion from surfaces.
Since the presence of torsion on the target manifold \(N\) does not affect the conformal structure on the domain \(M\),
Dirac-harmonic maps with torsion share many nice properties with ``usual'' Dirac-harmonic maps.
\begin{Lem}
\label{conformal-invariance}
In two dimensions the functional \(E_{\scriptscriptstyle{Tor}}(\phi,\psi)\) is conformally invariant.
\end{Lem}
\begin{proof}
It is well-known that the following terms are invariant under conformal transformations
\[
\int_M|d\phi|^2,\qquad \int_M\langle\psi,\D\psi\rangle,\qquad \int_M|\psi|^4
\]
and thus the energy functional \(E_{\scriptscriptstyle{Tor}}(\phi,\psi)\) is conformally invariant.
For more details, the reader may take a look at Lemma 3.1 in \cite{MR2262709}.
\end{proof}
For both harmonic and Dirac-harmonic map there exists a quadratic holomorphic differential,
we can find something similar here. Thus, let \((\phi,\psi)\) be a Dirac-harmonic map with torsion.
On a small domain \(\tilde{M}\) of \(M\) we choose a local isothermal
parameter \(z=x+iy\) and set
\begin{align}
\label{hopf-differential}
T(z)dz^2=&(|\phi_x|^2-|\phi_y|^2-2i\langle\phi_x,\phi_y\rangle \\
\nonumber&+\langle\psi,\partial_x\cdot\tilde{\nabla}^{Tor}_{\partial_x}\psi\rangle-i\langle\psi,\partial_x\cdot\tilde{\nabla}^{Tor}_{\partial_y}\psi\rangle)dz^2
\end{align}
with \(\partial_x=\frac{\partial}{\partial x}\) and \(\partial_y=\frac{\partial}{\partial y}\).

By varying \(E_{\scriptscriptstyle{Tor}}(\phi,\psi)\) with respect to the metric \(h_{\alpha\beta}\) of the domain \(M\),
we obtain the \emph{energy-momentum tensor}:
\begin{equation}
T_{\alpha\beta}=2\langle d\phi(e_\alpha),d\phi(e_\beta)\rangle-\delta_{\alpha\beta}|d\phi|^2
+\langle\psi,e_\alpha\cdot\nablattor_{e_\beta}\psi\rangle.
\end{equation}
It is clear that \(T_{\alpha\beta}\) is symmetric and traceless, when \((\phi,\psi)\)
is a Dirac-harmonic map with torsion.

\begin{Prop}
Let \((\phi,\psi)\) be a Dirac-harmonic map with torsion. Then the energy momentum tensor
is covariantly conserved, that is
\[
\nabla_{e_\alpha}T_{\alpha\beta}=0 .
\]
\end{Prop}
\begin{proof}
We choose a local orthonormal basis of \(TM\) with \(\nabla_{e_\alpha}e_\beta=0\) at the considered point.
By a direct calculation, using the skew-adjointness of the endomorphism \(A\), we obtain
\begin{align*}
\nabla_{e_\alpha}(2\langle d\phi(e_\alpha),d\phi(e_\beta)\rangle-\delta_{\alpha\beta}|d\phi|^2)=&2\langle\tau(\phi),d\phi(e_\beta)\rangle \\
=&2\langle\cR(\phi,\psi),d\phi(e_\beta)\rangle \\
&+2\langle F^{\scriptscriptstyle Tor}(\phi,\psi),d\phi(e_\beta)\rangle.
\end{align*}
Again, calculating directly, we get
\begin{align*}
\nabla_{e_\alpha}\langle\psi,e_\alpha\cdot\nablattor_{e_\beta}\psi\rangle=
&\langle\tilde{\nabla}^{\scriptscriptstyle LC}_{e_\alpha}\psi,e_\alpha\cdot\nablattor_{e_\beta}\psi\rangle +\langle\psi,\D(\nablattor_{e_\beta}\psi)\rangle \\
=&\langle A(d\phi(e_\alpha),e_\alpha\cdot\psi)),\nablattor_{e_\beta}\psi\rangle +\langle\psi,\D(\nablattor_{e_\beta}\psi)\rangle\\
=&\langle\psi,\Dtor(\nablattor_{e_\beta}\psi)\rangle,
\end{align*}
where we used that \(\psi\) is a solution of \eqref{euler-lagrange-psi-torsion}.
On the other hand, we find
\begin{align*}
\langle\psi,\Dtor\nablattor_{e_\beta}\psi\rangle=&\langle\psi,\nablattor_{e_\beta}\underbrace{\Dtor\psi}_{=0}\rangle
+\underbrace{\langle\psi,e_\alpha\cdot R^{\Sigma M}(e_\alpha,e_\beta)\psi\rangle}_{=\frac{1}{2}\langle\psi,\operatorname{Ric}(e_\beta)\cdot\psi\rangle=0} \\
&+\langle\psi,e_\alpha\cdot R^{N}_{\scriptscriptstyle Tor}(d\phi(e_\alpha),d\phi(e_\beta))\psi\rangle \\
=&-2\langle\cR(\phi,\psi),d\phi(e_\beta)\rangle-2\langle F^{\scriptscriptstyle Tor}(\phi,\psi),d\phi(e_\beta)\rangle.
\end{align*}
Adding up the different contributions then yields the assertion.
\end{proof}

\begin{Prop}
\label{prop-energy-momentum-preserved}
The quadratic differential \(T(z)dz^2\) is holomorphic.
\end{Prop}
\begin{proof}
This follows directly from the last Lemma.
\end{proof}

\begin{Lem}
The square of the twisted Dirac operator \(\Dtor\) satisfies the following Weitzenböck formula
\begin{align}
(\Dtor)^2\psi=&-\tilde{\Delta}^{\scriptscriptstyle Tor}\psi+\frac{R}{4}\psi+\frac{1}{2}e_\alpha\cdot e_\beta\cdot R^N(d\phi(e_\alpha),d\phi(e_\beta))\psi \\
\nonumber & +e_\alpha\cdot e_\beta\cdot\big((\nabla_{d\phi(e_\alpha)}A)(d\phi(e_\beta),\psi)+A(d\phi(e_\alpha),A(d\phi(e_\beta),\psi))\big),
\end{align}
where \(\tilde{\Delta}^{\scriptscriptstyle Tor}\) denotes the connection Laplacian on \(\Sigma M\otimes\phi^{-1}TN\).
\end{Lem}
\begin{proof}
This follows from a direct calculation or from the general Weitzenböck formula for twisted Dirac operators,
see for example \cite{MR1031992}, p.\ 164, Theorem 8.17 and \eqref{curvature-tensors}.
\end{proof}
As a next step we rewrite the Euler-Lagrange equations, for more details see \cite{MR2506243}.
By the Nash embedding theorem we can embed \(N\) isometrically in some \(\R^q\) of sufficient high dimension \(q\).
We then have that \(\phi\colon M\to\R^q\) with \(\phi(x)\in N\). The vector spinor \(\psi\)  becomes
a vector of untwisted spinors \(\psi^1,\psi^2,\ldots,\psi^q\), more precisely \(\psi\in\Gamma(\Sigma M\otimes T\R^q)\).
The condition that \(\psi\) is along the map \(\phi\) is now encoded as
\[
\sum_{i=1}^q\nu^i\psi^i=0\qquad \text{for any normal vector }\nu \text{ at } \phi(x).
\]
If we think of the torsion tensor \(A(\cdot,\cdot)\) as an endomorphism on \(TN\) we
can extend it to the ambient space \(\R^q\) by parallel transport.

\begin{Lem}
Assume that \(N\subset\R^q\). Moreover, assume that \(\phi\colon M\to\R^q\) and \(\psi\colon M\to\Sigma M\otimes T\R^q\).
Then the Euler-Lagrange equations acquire the form
\begin{align}
\label{phirq}-\Delta\phi=&\sff(d\phi,d\phi)+P(\sff(e_\alpha\cdot\psi,d\phi(e_\alpha)),\psi)+F^{\scriptscriptstyle Tor}(\phi,\psi),\\
\label{psirq}\p\psi=&\sff(d\phi(e_\alpha),e_\alpha\cdot\psi)+A(d\phi(e_\alpha),e_\alpha\cdot\psi),
\end{align}
where \(\sff\) denotes the second fundamental form in \(\R^q\) and \(P\) the shape operator.
\end{Lem}

\section{Analytic Aspects of Dirac-harmonic maps with torsion}
In this section we study analytic aspects of Dirac-harmonic maps with torsion.
This includes the regularity of solutions as well as the removal of isolated singularities.

\subsection{Regularity of solutions}
First of all, we need the notion of a weak solution of (\ref{euler-lagrange-phi-torsion}) and (\ref{euler-lagrange-psi-torsion}). 
Therefore, we define
\begin{align*}
\chi(M,N):=&\{(\phi,\psi)\in W^{1,2}(M,N)\times W^{1,\frac{4}{3}}(M,\Sigma M\otimes\phi^{-1}TN) \\
&\hspace{0.5cm}\text{ with } (\ref{phirq}) \text{ and } (\ref{psirq}) \text{ a.e.}\}.
\end{align*}
\begin{Dfn}[Weak Dirac-harmonic Map with torsion]
A pair \((\phi,\psi)\in\chi(M,N)\) is called \emph{weak Dirac-harmonic map with torsion} from \(M\) to \(N\) if and only if the pair \((\phi,\psi)\) 
solves (\ref{phirq}) and (\ref{psirq}) in a distributional sense.
\end{Dfn}
Note that the analytic structure of Dirac-harmonic maps with torsion is the same as the one of Dirac-harmonic maps
\begin{align*}
-\Delta\phi &\leq C(|d\phi|^2+|d\phi||\psi|^2), \\
\p\psi &\leq C|\psi||d\phi|. 
\end{align*}
Thus, the regularity theory developed for Dirac-harmonic maps can easily be applied.
More precisely, we may use the following (where \(D\) denotes the unit disc)

\begin{Satz}
Let \((\phi,\psi)\colon D\to N\) be a weak Dirac-harmonic map with torsion.
If \(\phi\) is continuous, then the pair \((\phi,\psi)\) is smooth.
\end{Satz}
This was proved in \cite{MR2176464}, Theorem 2.3, for Dirac-harmonic maps and can easily
be generalized to our case. Hence, we have to ensure the continuity of the map \(\phi\).
Thus, we will apply the following result due to Rivi\`ere (see \cite{MR2285745}):
\begin{Satz}
\label{theorem-riviere}
For every \(B={B_{~j}^i},1\leq i,j\leq q\)
in \(L^2(D,so(q)\otimes\R^2)\)(that is for all \(i,j\in 1,\ldots q, B_{~j}^i\in L^2(D,\R^2)\)
and \(B_{~j}^i=-B_{~i}^j\)), every \(\phi\in W^{1,2}(D,\R^q)\) solving
\begin{equation}
-\Delta\phi=B\cdot\nabla\phi
\end{equation}
is continuous. 
The notation should be understood as \(-\Delta \phi^i=\sum_{j=1}^qB_{~j}^i\cdot\nabla\phi^j\) for all \(1\leq i\leq q\).
\end{Satz}
To apply Theorem \ref{theorem-riviere} we further rewrite the Euler-Lagrange equations.
We will follow the presentation in \cite{MR2506243} for Dirac-harmonic maps.
We denote coordinates in the ambient space \(\R^q\) by \((y^1,y^2,\ldots,y^{q})\).
Let \(\nu_l,l=n+1,\ldots,q\) be an orthonormal frame field for the normal bundle \(T^\perp N\).
In addition, let \(D\) be a domain in \(M\) and consider a weak Dirac-harmonic map with torsion \((\phi,\psi)\in\chi(M,N)\).
We choose local isothermal coordinates \(z=x+iy\), set \(e_1=\partial_x\), \(e_2=\partial_y\)
and use the notation \(\phi_\alpha=d\phi(e_\alpha)\). 
The term involving the second fundamental form can be rewritten as
\begin{equation}
\label{skew-sff}
\sff^m(\phi_\alpha,\phi_\alpha)=\phi^i_\alpha\phi^j_\alpha\left(\frac{\partial\nu_l^i}{\partial y^j}\nu_l^m-\frac{\partial\nu_l^m}{\partial y^j}\nu_l^i\right),\qquad m=1,2,\ldots,q,
\end{equation}
see for example \cite{MR2285745}.
Following \cite{BVP}, p.7, the term on the right hand side of (\ref{phirq}) involving the shape operator can also
be written in a skew-symmetric way, namely
\begin{align}
\label{skew-p}
\operatorname{Re} P^m(\sff(\phi_\alpha,e_\alpha\cdot\psi),\psi)=& \\
\nonumber\phi^i_\alpha\langle\psi^k,e_\alpha\cdot\psi^j\rangle &
\Bigg(\bigg(\frac{\partial\nu_l}{\partial y^j}\bigg)^{\top,i}\bigg(\frac{\partial\nu_l}{\partial y^k}\bigg)^{\top,m}-\bigg(\frac{\partial\nu_l}{\partial y^k}\bigg)^{\top,i}\bigg(\frac{\partial\nu_l}{\partial y^j}\bigg)^{\top,m}\Bigg).
\end{align}
Here, \(\top\) denotes the projection map \(\top\colon\R^q\to T_yN\).

After these preparations we may now state the following
\begin{Prop}
Let \((M,h)\) be a closed Riemannian spin surface and let \(N\) be a compact Riemannian manifold. Assume that \((\phi,\psi)\in\chi(M,N)\) is a weak solution of (\ref{phirq})
and (\ref{psirq}). Let \(D\) be a simply connected domain of \(M\). Then there exists \(B^m_{~i	}\in L^2(D,so(q)\otimes\R^2)\)
such that
\begin{equation}
-\Delta\phi^m=B^m_{~i}\cdot\nabla\phi^i
\end{equation}
holds.
\end{Prop}
\begin{proof}
By assumption \(N\subset\R^{q}\) is compact, we denote its unit normal field by \(\nu_l,l=n+1,\ldots,q\).
Exploiting the skew-symmetry of (\ref{skew-sff}), (\ref{skew-p}), we denote
\[
B^m_{~i}=\begin{pmatrix}
f^m_{~i} \\
g^m_{~i}
\end{pmatrix},
\qquad
i,m=1,2,\ldots,q
\]
with
\begin{align*}
f^m_{~i}:=&\big(\frac{\partial\nu_l^i}{\partial y^j}\nu_l^m-\frac{\partial\nu_l^m}{\partial y^j}\nu_l^i\big)\phi^j_x\\
&+\langle\psi^k,\partial_x\cdot\psi^j\rangle_{\Sigma M}
\bigg(\bigg(\frac{\partial\nu_l}{\partial y^j}\bigg)^{\top,i}\bigg(\frac{\partial\nu_l}{\partial y^k}\bigg)^{\top,m}-\bigg(\frac{\partial\nu_l}{\partial y^k}\bigg)^{\top,i}\bigg(\frac{\partial\nu_l}{\partial y^j}\bigg)^{\top,m}\bigg)\\
&+\frac{1}{2}\langle\psi^k,\partial_x\cdot\psi^j\rangle_{\Sigma M}(\nabla^mA_{ijk}-\nabla_iA^m_{~~~jk}+A^m_{~~~~rk}A_{ij}^{\hspace{0.25cm}r}-A_{irk}A^{m\hspace{0.1cm}r}_{~~~~j})
\end{align*}
and we get the same expression for \(g^m_{~i}\) with \(x\) changed to \(y\).
Thus, we can write \eqref{phirq} in the following form
\[
-\Delta\phi^m=B^m_{~i}\cdot\nabla\phi^i.
\]
It remains to show that \(B^m_{~i}\in L^2(D,so(q)\otimes\R^2)\).
This follows directly since the pair \((\phi,\psi)\) is a weak solution of \eqref{phirq}, \eqref{psirq} 
and the Sobolev embedding \(|\psi|_ {L^4}\leq C|\psi|_{W^{1,\frac{4}{3}}}\).
The skew-symmetry of \(B^m_{~i}\) can be read of from its definition.
\end{proof}

\begin{Cor}
Let \((M,h)\) be a closed Riemannian spin surface and moreover, let \(N\subset \R^{q}\) be a compact manifold.
Suppose that \((\phi,\psi)\in\chi(M,N)\) is a weak Dirac-harmonic map with torsion. 
Then by the last Proposition and Theorem \ref{theorem-riviere}, 
we may deduce that \(\phi^m\) is continuous, \(m=1,2,\ldots,q\), hence \(\phi\in C^0(M,N)\).
\end{Cor}

\begin{Bem}
If we would consider a torsion contribution in the tension field, we could still deduce the continuity of the map \(\phi\)
due to the skew-adjointness of the endomorphism \(A\).
\end{Bem}
We may summarize our considerations by the following
\begin{Satz}
Let \((\phi,\psi)\colon D\to N\) be a weak Dirac-harmonic map with torsion.
Then the pair \((\phi,\psi)\) is smooth.
\end{Satz}
\subsection{Removable Singularity Theorem}
In this section we want to prove a removable singularity theorem for Dirac-harmonic maps
with torsion. More precisely, we want to show that solutions \((\phi,\psi)\) of (\ref{euler-lagrange-phi-torsion}) and (\ref{euler-lagrange-psi-torsion})
cannot have isolated singularities, whenever a certain energy is finite.
It is well known that such a theorem holds for both harmonic maps \cite{MR604040}
and also Dirac-harmonic maps \cite{MR2262709}. 
Let us define the following ``energy'':
\begin{Dfn}
Let \(U\) be a domain on \(M\). We define the energy of the pair \((\phi,\psi)\) on \(U\)
by
\begin{equation}
E(\phi,\psi,U):=\int_U(|d\phi|^2+|\psi|^4).
\end{equation}
\end{Dfn}
This energy is conformally invariant and thus plays an important role. 

First of all, we need some local energy estimates (with the unit disc \(D\)). 
\begin{Satz}
\label{theorem-energy-local}
Let \((M,h)\) be a closed Riemannian spin surface and \((N,g)\) a
compact Riemannian manifold. Assume that the pair \((\phi,\psi)\) is a Dirac-harmonic map with torsion.
There is a small constant \(\epsilon>0\) such that if the pair \((\phi,\psi)\) satisfies
\begin{equation}
\int_D(|d\phi|^2+|\psi|^4)<\epsilon,
\end{equation}
then
\begin{equation}
|d\phi|_{C^k(D_\frac{1}{2})}+|\psi|_{C^k(D_\frac{1}{2})} \leq C(|d\phi|_{L^2(D)}+|\psi|_{L^4(D)}),
\end{equation}
where the constant \(C\) depends on \(N\) and \(k\).
\end{Satz}
Since the presence of torsion does not affect the analytic structure of the Euler-Lagrange equations 
the same proof as for Theorem 4.3 in \cite{MR2262709} still holds.

The behaviour of \((\phi,\psi)\) near a singularity can be described by the following 
\begin{Cor}
\label{corollary-energy-local}
There is an \(\epsilon>0\) small enough such that if the pair \((\phi,\psi)\) is a smooth solution of (\ref{euler-lagrange-phi-torsion})
and (\ref{euler-lagrange-psi-torsion}) on \(D\setminus\{0\}\) with finite energy \(E(\phi,\psi,D)<\epsilon\), 
then for any \(x\in D_{\frac{1}{2}}\) we have
\begin{align}
|d\phi(x)||x|\leq& C|d\phi|_{L^2(D_{2|x|})}, \\
|\psi(x)|^\frac{1}{2}|x|^\frac{1}{2}+|\nabla\psi(x)||x|^\frac{3}{2}\leq& |\psi|_{L^4(D_{2|x|})}.
\end{align}
\end{Cor}
\begin{proof}
Fix any \(x_0\in D\setminus\{0\}\) and define \((\tilde{\phi},\tilde{\psi})\) by
\[
\tilde{\phi}(x):=\phi(x_0+|x_0|x) \textrm{ and } \tilde{\psi}(x):=|x_0|^\frac{1}{2}\psi(x_0+|x_0|x).
\]
It is easy to see that \((\tilde{\phi},\tilde{\psi})\) is a smooth solution of (\ref{euler-lagrange-phi-torsion}) and (\ref{euler-lagrange-psi-torsion}) on \(D\)
with \(E(\tilde{\phi},\tilde{\psi},D)<\epsilon\). By application of Theorem \ref{theorem-energy-local}, we have
\[
|d\tilde{\phi}|_{L^\infty(D_\frac{1}{2})}\leq C|d\tilde{\phi}|_{L^2(D)}, \qquad |\tilde{\psi}|_{C^1(D_\frac{1}{2})}\leq C|\tilde{\psi}|_{L^4(D)}
\]
and scaling back yields the assertion.
\end{proof}

\begin{Prop}
\label{prop-polar-coordinates}
Let \((\phi,\psi)\) be a smooth Dirac-harmonic map with torsion on \(D\setminus\{0\}\)
satisfying \(E(\phi,\psi,D)<\epsilon\). Then we have
\begin{align*}
\int_0^{2\pi}\frac{1}{r^2}\big|\frac{\partial\phi}{\partial\theta}\big|^2d\theta&=\int_0^{2\pi}\big|\frac{\partial\phi}{\partial r}\big|^2+\langle\psi,\partial_r\cdot\frac{\tilde{\nabla}\psi}{\partial r}\rangle
+\langle\psi,A(\frac{\partial\phi}{\partial r},\partial_r\cdot\psi)\rangle\\
&=\int_0^{2\pi}\big|\frac{\partial\phi}{\partial r}\big|^2-\frac{1}{r^2}\langle\psi,\partial_\theta\cdot\frac{\tilde{\nabla}\psi}{\partial\theta}\rangle
-\frac{1}{r^2} \langle\psi,A(\frac{\partial\phi}{\partial\theta},\partial_\theta\cdot\psi)\rangle,
\end{align*}
where \((r,\theta)\) are polar coordinates on the disc \(D\) centered around the origin.
\end{Prop}
 
\begin{proof}
By Proposition (\ref{prop-energy-momentum-preserved}) we know that
\begin{align*}
T=&|\phi_x|^2-|\phi_y|^2-2i\langle\phi_x,\phi_y\rangle
+\langle\psi,\partial_x\cdot\nablattor_{\partial_x}\psi\rangle-i\langle\psi,\partial_x\cdot\nablattor_{\partial_y}\psi\rangle 
\end{align*}
is holomorphic on \(D\setminus\{0\}\) with \(z=x+iy\in D\). By the last Corollary we also know that
\begin{align*}
|\psi||\nablattor\psi|\leq C(|\psi||\nabla^{\scriptscriptstyle\Sigma M}\psi|+|d\phi||\psi|^2)\leq C|z|^{-2},
\qquad |d\phi|^2\leq C|z|^{-2},
\end{align*}
which, altogether, gives \(|T(z)|\leq Cz^{-2}\).
Moreover, it is easy to see that \(\int_D|T(z)|<\infty\).
Hence, we may follow that \(zT(z)\) is holomorphic on the disc \(D\) and by Cauchy's integral theorem we deduce
\begin{equation}
0=\operatorname{Im}\int_{|z|=r}zT(z)dz=\int_0^{2\pi}\operatorname{Re}(z^2T(z))d\theta.
\end{equation}
Moreover, a direct calculation gives
\begin{align}
\operatorname{Re}(z^2T(z))=r^2\big|\frac{\partial\phi}{\partial r}\big|^2-\big|\frac{\partial\phi}{\partial\theta}\big|^2-\langle\psi,\partial_\theta\cdot\frac{\nablattor\psi}{\partial\theta}\rangle,
\end{align}
which finally proves the result.
\end{proof}

Now we are in the position to state the
\begin{Satz}[Removable Singularity Theorem]
Let \((\phi,\psi)\) be a solution of (\ref{euler-lagrange-phi-torsion}) and (\ref{euler-lagrange-psi-torsion}), which is smooth on \(U\setminus\{p\}\) for some \(p\in U\).
If \((\phi,\psi)\) has finite energy \(E(\phi,\psi,D)\), then \((\phi,\psi)\) extends to a smooth solution on \(U\). 
\end{Satz}

\begin{proof}
With the help of the last Proposition the same proof as for Theorem 4.6 in \cite{MR2262709} can be applied.
\end{proof}

\section{Dirac harmonic with curvature term and torsion}
As we have seen in the introduction the full \((1,1)\) non-linear supersymmetric sigma model 
studied in quantum field theory involves an additional curvature term in the energy functional,
namely
\begin{equation}
\label{energy-curvature}
E_c(\phi,\psi)=\frac{1}{2}\int_M|d\phi|^2+\langle\psi,\D\psi\rangle+\frac{1}{6}\langle R^N(\psi,\psi)\psi,\psi\rangle.
\end{equation}
Here, the indices are contracted as follows
\[
\langle R^N(\psi,\psi)\psi,\psi\rangle = R_{ijkl}\langle\psi^i,\psi^k\rangle\langle\psi^j,\psi^l\rangle,
\]
which ensures that the action is real-valued.
The factor \(\frac{1}{6}\) in front of the curvature term is required by supersymmetry, see \cite{MR1707282}, p.78.
This functional has already been analyzed in \cite{MR2370260}.

\begin{Prop}[\cite{MR2370260}]
The critical points of the energy functional (\ref{energy-curvature}) are given by
\begin{align}
\label{phi-curvature}\tau(\phi)=&\cR(\phi,\psi)+\tilde{\cR}(\psi), \\
\label{psi-curvature}\D\psi=&\frac{1}{3}R^N(\psi,\psi)\psi
\end{align}
with the curvature terms
\begin{equation*}
\tilde{\cR}(\psi)=\frac{1}{12}\langle(\nabla R)^\sharp(\psi,\psi)\psi,\psi\rangle,\qquad
\cR(\phi,\psi)=\frac{1}{2}R^N(e_\alpha\cdot\psi,\psi)d\phi(e_\alpha).
\end{equation*}
Here, \(\sharp\colon T^\ast N\to TN\) denotes the musical isomorphism.
\end{Prop}
Solutions \((\phi,\psi)\) of the system \eqref{phi-curvature}, \eqref{psi-curvature} are called \emph{Dirac-harmonic maps with curvature term}.
\par\medskip
At present very little is known about the properties of Dirac-harmonic maps with curvature term.
Compared to Dirac-harmonic maps the main difference arises in the fact that Dirac-harmonic maps
with curvature term constitute a coupled system of two non-linear equations.
This makes the analysis of solutions of \eqref{phi-curvature} and \eqref{psi-curvature}
substantially harder. However, \eqref{psi-curvature} has an interesting limit.
In the case that the map \(\phi\) is trivial \eqref{psi-curvature} gives rise
to the spinorial Weierstrass representation of surfaces. The analytic aspects of this
equation are investigated in \cite{MR2390834} and \cite{MR2661574}.
\par\medskip
The functional \eqref{energy-curvature} also has a natural extension to target spaces with torsion
(see again \(E_{ASYM}(\Phi)\) in the introduction).
In this case, we have to specify \(H(\psi,\psi,\psi,\psi)\) in \eqref{easym}.
If we replace the curvature tensor in \eqref{energy-curvature} with the curvature tensor
of a connection with torsion and contract the indices the same way, then in general the action will no longer be real-valued.
This is due to the fact that the curvature tensor is not symmetric under swapping the first with the second pair of indices,
see Rem. 2.3. in \cite{MR2322400}. However, if we stick to totally anti-symmetric torsion
and impose the additional condition that the torsion is parallel, then the curvature tensor
has the necessary symmetries to obtain a real-valued action.
When studying the critical points of this functional one obtains a set of equations that has
the same analytic structure as Dirac-harmonic maps with curvature term.
  
\begin{Bem}
Together with the analysis performed in \cite{magnetic}
one gets a full description of the full \((1,1)\) non-linear supersymmetric sigma model from
the perspective of differential geometry.
\end{Bem}

\bibliographystyle{alpha}
\bibliography{mybib}
\end{document}